\documentclass[12pt]{article}
 
\input xy
\xyoption{all}
 
\usepackage{amssymb}
\usepackage{amsmath}
\usepackage{amscd}
\usepackage{latexsym}

\begin{document}

\newtheorem{lemma}{Lemma}[section]
\newtheorem{prop}[lemma]{Proposition}
\newtheorem{cor}[lemma]{Corollary}
 
\newtheorem{thm}[lemma]{Theorem}
\newtheorem{con}{Conjecture}

\newtheorem{rem}[lemma]{Remark}
\newtheorem{rems}[lemma]{Remarks}
\newtheorem{defi}[lemma]{Definition}
 \newtheorem{ex}[lemma]{Example}

\newcommand{\C}{\mathbb C}
\newcommand{\R}{\mathbb R}
\newcommand{\Q}{\mathbb Q}
                                                                                
\newcommand{\Z}{\mathbb Z}
\newcommand{\N}{\mathbb N}
\newcommand{\MaxG}{C^{\ast}G}

\newcommand{\G}{\mathcal G}                                                                  

\newcommand{\vNG}{{\mathcal N}\!(G)}
\newcommand{\vNH}{{\mathcal N}\!(H)}
\newcommand{\ldg}{\ell^2(G)}
                                                                                
\newcommand{\ldh}{\ell^2(H)}
\newcommand{\cg}{\C{G}}
\newcommand{\zg}{\Z{G}}
\newcommand{\qg}{\Q{G}}

\newcommand{\lte}{\chi^{(2)}}
                                                                                
\newcommand{\Tr}{\rm Tr}
\newcommand{\Hn}{{\mathcal H}^{(2)}_n}
\newcommand{\tr}{\operatorname{tr}_{\mathcal N}}
                                                                                
\newcommand{\trg}{\operatorname{tr}_{{\mathcal N}G}}
\newcommand{\trh}{\operatorname{tr}_{{\mathcal N}H}}
                                                                                
\newcommand{\dvNG}{\operatorname{dim}_{G}}
\newcommand{\HS}{\operatorname{HS}}

\title{Hattori-Stallings Trace and\\ Euler Characteristics for Groups
\footnote{Revised December 16, 2004}}

                                                                                
\author{Indira Chatterji\footnote{Partially supported by the Swiss National Science Foundation} and Guido Mislin}
\date{{}}
\maketitle
\section*{Introduction}
For $G$ a group and $P$ a finitely generated projective
module over the integral group ring, 
Bass conjectured in \cite{B} that the
Hattori-Stallings rank of $P$
should vanish on elements different from $1\in G$,
and
proved it in many cases such as torsion-free linear groups. Later, this
conjecture has been proved for many more groups, notably by Eckmann \cite{E},
Emmanouil \cite{Emmanouil} and Linnell \cite{Linnell}. The latest 
advances are given in \cite{BCM} and
the first section of the present paper is a quick survey of
the Bass conjecture, together with an outline of the proof of the
main result of \cite{BCM}. 

In most cases, one proves a stronger conjecture,
which asserts that the Hattori-Stallings rank of a finitely 
generated projective module over the
complex group ring should vanish on elements of infinite order
(that this conjecture is indeed stronger follows from
Linnell's work \cite{Linnell}). Given a
group $G$ of type FP over $\C$, its complete Euler characteristic 
$E(G)$ is the Hattori-Stallings rank
of an alternating sum of finitely generated projective modules over $\cg$, and on the elements of
finite order, one could then ask of what the values do depend.
It is Brown in \cite{Brown}
who first studied that question, proving a formula in many cases. In Section \ref{euler} we shall explain 
the basics to understand Brown's formula and propose 
Conjecture 1 below as a generalization. Our generalization amounts to 
putting Brown's work in the context of $L^2$-homology (not available at the time where \cite{Brown} has been written),
and applies to cases where Brown's
formula is not available. 
%
\begin{con}\label{conjecture}Let $G$ be a group of type FP over $\C$ such that the centralizer of every
element of finite order in $G$ has finite $L^2$-Betti numbers. Then for
every $s\in G$
\begin{equation}\label{main}E(G)(s)=\lte(C_G(s)),
\end{equation}
where $E(G)(s)$ is the $s$-component of the complete Euler characteristic of $G$ and $\lte(C_G(s))$ is the $L^2$-Euler characteristic of the centralizer of $s$ in $G$.
\end{con}
This formula, as opposed to the Bass conjecture, has nice stability 
properties that we
discuss in Section \ref{stability}. We describe in Section \ref{Binfini} 
a class of groups containing all $G$ with cocompact
$\underline{E}G$ for which Formula (\ref{main}) 
holds.
It is straightforward (see Lemma \ref{chi=e} below) 
that Formula (\ref{main}) always
holds for $s=1$.
%
If we write $\chi(G)$ for the naive Euler characteristic
$\sum(-1)^i \operatorname{dim}_\C H_i(G;\C)$ then for $G$
satisfying Conjecture \ref{conjecture}, we find (cf.~Corollary \ref{chis})
\begin{equation}\label{genAtiyah}
\chi(G)=\sum_{[s]\in [G]}\lte(C_G(s)).
\end{equation}
If $K(G,1)$ is a finite complex, then $G$ satisfies Conjecture
\ref{conjecture} and $G$ is necessarily torsion-free
so that Formula (\ref{genAtiyah}) reduces to Atiyah's celebrated theorem:
$\chi(G)=\lte(G)$.
\section{Review of the Bass conjecture}\label{ReviewBass}
For a group $G$ we denote by $\HS: K_0(\C G)\to HH_0(\C G)=\bigoplus_{[G]}\C$
the {\sl Hattori-Stallings Trace}; $[G]$ stands for the set of conjugacy classes of $G$. If $P$ denotes a finitely generated projective $\cg$-module and $[P]\in K_0(\cg)$ the corresponding element, we write
$$\HS(P):=\HS([P])=\sum_{[s]\in [G]}\HS(P)(s)\cdot [s]\in\bigoplus_{[G]}\C$$
with $\HS(P)(s)$ depending on the conjugacy class $[s]$ of $s\in G$ only. Therefore, we can think of $\HS(P): G\to \C$ as a class function. 
It is well-known that for $s\in G$ a {\sl central} element of 
infinite order, one has $\HS(P)(s)=0$. More generally, 
if $C_G(s)$ denotes the centralizer of $s\in G$, and $[G:C_G(s)]$ 
is finite and $s$ has infinite order, then $\HS(P)(s)=0$, 
because $\HS(P|_{C_G(s)})(s)=\HS(P)(s)$ (in general, if $H<G$
is a subgroup of finite index and $s\in H$, then $\HS(P|_H)(s)$ =
$[C_G(s):C_H(s)]\HS(P)(s)$, see \cite[Corollary 6.3]{B} or 
Chiswell's notes \cite{chiswell}). 

Another very useful result in this context goes back to 
Bass (cf.~\cite{B}, Proposition 9.2), and states that if $\HS(P)(s)\ne 0$ 
then there is an $N>0$ such that for almost all primes $p$, 
the elements $s^{p^N}$ are conjugate to $s$; note that in 
case $s$ has infinite order, this implies that for almost 
all primes $p$, $s$ is contained in a subgroup of $G$ which 
is isomorphic to $\Z[1/p]$.
According to Bass in \cite{B}, the following general 
vanishing theorem ought to be true:
\begin{con}[Bass Conjecture over $\C$]\label{BC}
For $P$ a finitely generat\-ed projective $\cg$-module and $s\in G$ an element of infinite order, $\HS(P)(s)= 0$.\end{con}
The Bass Conjecture over $\C$ is known to hold for many groups, including:
\begin{itemize}
\item linear groups (cf.~Bass \cite{B}), which includes Braid Groups (they are linear: \cite{K} and \cite{Bigelow})
\item groups with $\operatorname{cd}_\C \le 2$ (cf.~Eckmann \cite{E}; see also Emmanouil \cite{Emmanouil}, as well as \cite{EP} 
for more general results using techniques of cyclic homology), which includes one-relator groups and knot groups
\item subgroups of semihyperbolic groups (this follows from results of Alonso
and Bridson \cite{AB}, see also \cite{BE} or the discussion in \cite{MV}, 
following Corollary 7.17 of Part 1; for the definition of 
semihyperbolic groups the reader is referred to \cite{BH}); 
these include (subgroups of) word hyperbolic groups and 
cocompact CAT(0)-groups
\item mapping class groups $\Gamma_g$ of closed surfaces of 
genus $g$ (cf.~Corollary 7.17 (Part 1) of \cite{MV})
\item amenable groups, more generally groups satisfying the Bost Conjecture (for the Bost Conjecture, see \cite{Lafforgue} and \cite{Skandalis}) have been shown to satisfy the Bass Conjecture over $\C$ in \cite{BCM}. For instance groups which have the {\sl Haagerup Property} (also called {\sl a-T-menable groups}). We recall that a group is said to have the Haage\-rup Property if it admits an isometric, metrically proper affine action on some Hilbert space (for a discussion of such groups, see \cite{CCJJV}); the class of groups having the Haagerup property contains all countable groups which are extensions of amenable groups with free kernel, and is closed under subgroups, finite products, passing to the fundamental group of a countable, locally finite graph of groups with finite edge stabilizers (vertex stabilizers are assumed to have the Haagerup property), countable increasing unions, amalgamations $A\ast_B C$ with $A$ and $C$ both countable amen\-able and $B$ central in $A$ and $C$ (use Propositions 4.2.12 and 6.2.3 of \cite{CCJJV}) and passing to finite index supergroups. Groups which act metrically properly and isometrically on a uniformly locally finite, weakly $\delta$-geodesic and strongly $\delta$-bolic space (see \cite{KS} and \cite{Lafforgue}); examples of groups satisfying these conditions are word hyperbolic groups (see \cite{MiYu}) and cocompact CAT(0)-groups.
\end{itemize}
A prominent class of groups for which Conjecture \ref{BC} is not known
 in general, is the class of profinite groups. However, if $G$ is any 
group and $Q$ a finitely generated projective $\zg$-module, then 
according to 
Linnell \cite{Linnell}, if $s\ne 1$ is such that 
$\HS(\C G\otimes_{\Z G}Q)(s)\ne 0$, then $s$ is contained in a subgroup 
of $G$ isomorphic to the additive group $\Q$ of rationals. This 
in particular implies that for $G$ profinite one has 
$\HS(\C G\otimes_{\Z G}Q)(s)=0$ for all $s\in G\setminus\{1\}$, 
because $\Q$ cannot be a subgroup of a profinite group.

We give an outline of the strategy for proving the main result of \cite{BCM}, which states that the Bost Conjecture implies the Bass Conjecture over $\C$ (see Theorem \ref{BOSTimpliesBASS} below). The Bost Conjecture asserts that the {\sl Bost assembly map}
$$\beta_0^G: K_0^G( \underline{E}G)\to K_0( \ell^1G)$$
is an isomorphism (see \cite{Lafforgue} and \cite{Skandalis}).
Here, the left hand side denotes the equivariant 
$K$-homology of the classifying space for proper actions of $G$, 
and the right hand side is the projective class group of the Banach 
algebra $\ell^1 G$ of summable complex valued functions on $G$.
\begin{thm}\label{BOSTimpliesBASS}
Suppose that $G$ satisfies the Bost Conjecture. Then $G$ satisfies the Bass Conjecture over $\C$.
\end{thm}
Before outlining the proof of Theorem \ref{BOSTimpliesBASS}, we need to address some auxiliary constructions. We can extend $\HS: K_0( \C G)\to \bigoplus_{[G]}\C$ to a trace $\HS^{1}: K_0(\ell^1 G)\to\prod_{[G]}\C$ as follows. If $[Q]\in K_0(\ell^1G)$, with $Q$ a finitely generated projective $\ell^1G$-module, we choose an idempotent $(n,n)$-matrix $M=(m_{ij})$ with entries in $\ell^1G$ 
representing $Q$ (i.e. $(\ell^1G)^n\cdot M\cong Q$ as 
left $\ell^1G$-modules), then we put
$$\HS^1(Q):= \HS^1([Q]):=\{\sum_{i=1}^n
\sum_{t\in[s]}m_{ii}(t)\}_{[s]\in[G]}\in\prod_{[G]}\C.$$
The $m_{ii}(t)$'s stand for the $t$-coefficients of $m_{ii}\in\ell^1G$, $1\le i\le n$. We will write $\HS^1(x)(s)$ for the $[s]$-component of $\HS^1(x),\, x\in K_0(\ell^1 G)$. One checks that $\HS^1$ is well-defined and compatible with the ususal Hattori-Stallings trace.
%
To get informations on $\HS^1$ via the Bost assembly map, we embed $G$ into an acyclic group of a very special kind. Recall that a group $G$ is called \emph{acyclic}, if $H_i(G;\Z)=0$ for $i>0$. As proved in \cite{BCM}, every group $G$ admits a functorial embedding into an acyclic group $A=A(G)$, which we call the {\sl pervasively acyclic hull of $G$}, satisfying the following:
\begin{itemize} 
\item For every finitely generated abelian subgroup $B<A$ the centralizer $C_{A}(B)$ is acyclic (such a group is called \emph{pervasively acyclic})
\item $A$ is countable if $G$ is and the induced map on conjugacy classes $[G]\to [A]$ is injective.
\end{itemize}
In this context, the important feature of a pervasively 
acyclic group $A$ is that its classifying space for proper 
actions is \emph{$K_0^{A}\otimes\Q$-discrete}, meaning that the
inclusion $\underline{E}A^0\to\underline{E}A$ of the
$0$-skeleton induces
a surjective map
$$K^{A}_0(\underline{E}A^0)\otimes\Q\to K_0^{A}(\underline{E}A)\otimes\Q,$$
see Corollary 3.9 of \cite{BCM}. 
In other words, all the information of the $A$-$CW$-complex $\underline{E}A$ captured by the equivariant $K$-homology is contained in its $0$-skeleton 
$\underline{E}A^0 = \coprod_{\alpha} A/A_\alpha,$
where $A_\alpha<A$ stands for a finite subgroup, corresponding to the 
stabilizer of some $0$-cell of $\underline{E}A$. 
The equivariant $K$-homology we use here is the one defined by 
Davis and L\"uck (cf.~\cite{DL}), arising from a spectrum over 
the orbit category of $G$. It is defined on the category of 
\emph{all} $G$-$CW$-complexes; on proper, cocompact $G$-$CW$-complexes, 
this {\sl representable} equivariant $K$-homology agrees with 
the one used in the original version of the Baum-Connes or Bost 
conjectures (see \cite{HP}) so that if $X$ is a proper, not necessarily cocompact 
$G$-$CW$-complex, then $K_*^G(X)=\operatorname{colim}_{Y\subset X, Y/G\, {\text{compact}}}K^G_*(Y),$
in accordance with the classical setup for the Baum-Connes and Bost conjectures. It follows that $K_0^G$ is fully additive so that
$$K_0^A(\underline{E}A^0)=\bigoplus_\alpha K_0^A(A/A_\alpha)=\bigoplus_\alpha K_0^{A_\alpha}(\{pt\})=\bigoplus_\alpha K_0(\ell^1A_\alpha)$$
and $K_0(\ell^1A_\alpha)\cong R_\C(A_\alpha)$, the additive group of the complex representation ring of the finite group $A_\alpha$. 

\medskip

{\sl Outline of the proof of Theorem \ref{BOSTimpliesBASS}}. Let $P$ be a finitely generated projective $\cg$-module and assume that $G$ satisfies the Bost Conjecture. Then $x:=[\ell^1G\otimes_{\cg} P]$ lies in the image of the Bost assembly map $\beta_0^G$ and $\HS^1(x)$ captures the information contained in $\HS(P)$. We embed $G$ into its pervasively acyclic hull $A(G)=:A$ and together with the standard embedding $\underline{E}A^0\to\underline{E}A$ 
of the $0$-skeleton this yields a commutative diagram
$$\begin{CD}
@. [P]\in K_0(\C G){\mbox{\phantom{xxxxx}}}@>{\HS}>>\bigoplus_{[G]}\C\\
@. @VVV @VVV\\
K_0^G(\underline{E}G)@>{\beta^G_0}>{\cong}> K_0(\ell^1G)@>{\operatorname{HS}^1}>>
\prod_{[G]}\C\\
@VVV @VVV @VVV\\
K_0^A(\underline{E}A)@>{\beta^A_0}>> K_0(\ell^1A)@>{\HS^1}>>\prod_{[A]}\C\\
@AAA @AAA @AAA\\
K_0^A(\underline{E}A^0)@>{\bigoplus\beta_0^{A_\alpha}}>>
\bigoplus_\alpha K_0(\ell^1A_\alpha)
@>{\bigoplus_\alpha\HS}>>\bigoplus_{\alpha,[A_\alpha]}\C.\\
\end{CD}
$$
Using the facts that $K_0^A(\underline{E}A^0)\to K_0^A(\underline{E}A)$ is 
rationally surjective (since $\underline{E}A$ is $K_0^A\otimes\Q$-discrete), 
and that the induced map $\prod_{[G]}\C \to\prod_{[A]}\C$ is injective, 
we conclude by diagram chasing
that $\HS^1(x)$ lies in the subspace of functions $[G]\to\C$, 
whose support is contained in the subset of those conjugacy 
classes of $G$, which are represented by elements of finite order. 
Therefore $\HS^1(x)(s)=0$ for $s\in G$ of infinite order, 
which implies that $\HS(P)(s)=0$ too, establishing the 
Bass conjecture over $\C$ for the group $G$.

\hfill QED
\section{Euler characteristics}\label{euler}
In this section we shall explain the basics to discuss Conjecture \ref{conjecture}. Let $G$ be a group of type FP over $\C$, meaning that there exists a resolution
$$P_*: 0\to P_n\to P_{n-1}\to\cdots\to P_0\to\C$$
with each $P_i$ finitely generated projective over $\C G$; in case
the $P_i$'s may be chosen to be finitely generated and free over $\C G$,
the $G$ is termed of type FF over $\C$. 
The element $W(G):=\sum_{i}(-1)^{i}[P_i]\in K_0(\C G)$ depends 
on $G$ only and we call it the \emph{Wall element}. 
Under the Hattori-Stallings trace, the Wall element $W(G)$ is mapped to 
$E(G)=\sum_{[s]\in [G]} E(G)(s) [s],$ 
the sum being taken over the set  $[G]$ of conjugacy classes $[s]$ 
of elements $s\in G$. This is the \emph{complete Euler 
characteristic of $G$} (see \cite{S}). If $G$ has a 
cocompact $\underline{E}G$, Conjecture \ref{conjecture} is 
true as it reduces to Brown's formula \cite{Brown} that we 
shall now discuss. For $G$ of type FP over $\C$, 
the \emph{Euler characteristic of $G$} (in the sense of Bass \cite{B} 
and Chiswell \cite{chiswell}) is given by
$e(G)=E(G)(1).$
Note also that $W(G)=0$ if and only if $e(G)=0$ and $G$ is of type FF
over $\C$.
Brown conjectures under suitable finiteness conditions for $G$ 
the following formula:
\begin{equation}\label{brown}
E(G)(s)=\left\{\begin{array}{cc}e(C_G(s)) & \hbox{if $s$ has finite order}\\
0 & \hbox{otherwise}\end{array}\right.
\end{equation}
and proves it in many cases, including groups with
cocompact $\underline{E}(G)$. Brown's assumptions always require 
in particular
$C_G(s)$ to be of
type FP over $\C$, and in this case we will show that Formula (\ref{main}) 
reduces to Brown's formula (\ref{brown}). To do this, we first recall the definition of $L^2$-Euler characteristic. For $i\in\N$, the $i$-th $L^2$-Betti 
number is defined as the von Neumann dimension of the $\vNG$-module $H_i\left(G;\vNG\right)$
$$\beta_i(G)=\dvNG H_i\left(G;\vNG\right)\in [0,\infty],$$
where $\vNG$ is the group von Neumann algebra of $G$ 
(see L\"uck's book \cite{L}). If $\sum(-1)^i\beta_i(G)$ converges, 
the \emph{$L^2$-Euler characteristic} is defined as
\begin{equation}\label{Euler}
\lte(G)=\sum_{i\in\N}(-1)^i\beta_i(G)\in\R.\end{equation}
In case $G$ is finite, $\lte(C_G(s))=1/|C_G(s)|$ and Formulae 
(\ref{main}) and (\ref{genAtiyah}) reduce to well-known results. 
With no finiteness restrictions imposed on $G$, 
one can find for any $r\in\R$ a group $G$ with $\lte(G)=r$. 
However, if $G$ is of type FP over $\C$ then $\lte(G)\in\Q$, 
as shown by the following.
\begin{lemma}{\label{chi=e}}
Suppose that a group $G$ is of type FP over $\C$. Then
$\lte(G)=e(G)$ and $e(G)$ is a rational number.\end{lemma}
{\it Proof.} Let $P_*: 0\to P_n\to P_{n-1}\to\cdots\to P_0\to\C$
be a projective resolution of type FP for $G$. Then
\begin{eqnarray*}\lte(G)&=&\sum_{i\in\N}(-1)^i\beta_i(G)=\sum_{i\in\N}(-1)^i\dvNG
\vNG\otimes_{\cg}P_i\\
&=& \sum_{i\in\N}(-1)^iHS(P_i)(1)=e(G).\end{eqnarray*}
Here we used the fact that for a finitely generated projective $\cg$-module $P$, 
$\dvNG\vNG\otimes_{\cg}P=HS(P)(1)$, which is actually 
just the Kaplansky trace of $P$; the Kaplansky trace of a finitely 
generated projective $\cg$-module is a rational number, 
by Zalesskii's theorem (see \cite{BV}).\hfill QED

\medbreak
The $L^2$-Betti numbers turn out to be often $0$. In particular we mention the following vanishing result.
\begin{thm}[Cheeger-Gromov, \cite{CG}]\label{CG}If $G$ contains an infinite normal amenable subgroup, then $\beta_i(G)=0$ for all $i\in\N$, and therefore $\lte(G)=0$.\end{thm}
This theorem immediately implies that for an arbitrary group $G$, the $L^2$-Euler characteristic of $C_G(s)$ is $0$ for all $s\in G$ of infinite order, so that one more evidence for Conjecture \ref{conjecture} is the following simple observation: \emph{If the Bass Conjecture over $\C$ holds for $G$, then Formula (\ref{main}) holds on elements of infinite order.}
Indeed, the Bass conjecture will say that the left hand side vanishes 
on elements of infinite order. The following fact on $L^{2}$-Euler characteristics will be used later, mainly for the case of subgroups $H<G$, with $G$ of type FP over $\C$. Since then $\operatorname{cd}_\C H<\infty$, 
the Euler characteristic $\lte(H)$ is well-defined 
if and only if all $L^2$-Betti numbers $\beta_i(H)$ are finite.
\begin{lemma}\label{prodchi}
Let $H$ and $K$ be groups with $\sum_i\beta_i(H)$ and $\sum_i\beta_i(K)$ convergent. Then $\lte(H\times K)=\lte(H)\lte(K)$.\end{lemma}
{\it Proof.} One uses the K\"unneth Formula for $L^{2}$-Betti numbers
\cite{CG}: $\beta_n(H\times K)=\sum_{i+j=n}\beta_i(H)\beta_j(K)$,
and takes the alternating sum; note that $\sum_n\beta_n(H\times K)$ is convergent so that $\lte(H\times K)$ is well-defined.\hfill QED

\medskip

A $1$-dimensional contractible $G$-$CW$-complex $T$ with vertex set $V$ and edge set $E$ (for short: a $G$-tree) is given by a $G$-push-out
$$\begin{CD}
(\coprod_{\beta\in E/G} G/G_\beta)\times S^0 @>>>
\coprod_{\alpha\in V/G} G/G_\alpha\\
@VVV @VVV\\
(\coprod_{\beta\in E/G} G/G_\beta)\times D^1@>>> T\\
\end{CD}$$
and the cellular chain complex of $T$ has the form
$$0\to\bigoplus_{\beta\in E/G}\C[G/G_\beta]\to\bigoplus_{\alpha\in V/G}\C [G/G_\alpha]\to\C.$$
The group $G$ is then the fundamental group of a graph of groups 
$\{G_\gamma\}_{\gamma\in I}$, $I=V/G\sqcup E/G$; the graph is called 
{\it finite}, if $I$ is a finite set (i.e.~if the action of 
$G$ on $T$ is cocompact). If $X$ is an arbitrary $G$-$CW$-complex,
we write $H_*(X;\vNG):=H_*(\vNG\otimes_{\Z G}C^{\text{cell}}_*(X))$ 
for its $L^2$-homology
so that $H_*(G;\vNG)=H_*(EG;\vNG)$.

\begin{lemma}\label{graph} 
Let $G$ be the fundamental group of a (not necessarily finite) graph of groups $\{G_\gamma\}_{\gamma\in I}$, where $I=V/G\sqcup E/G$. If for each of the groups $G_\gamma$ the series $\sum_i\beta_i(G_\gamma)$ is convergent and equals $0$ for almost all $\gamma\in I,$ then
$$\lte(G)=\sum_{\alpha\in V/G}\lte(G_\alpha) - \sum_{\beta\in E/G}\lte(G_\beta).$$
\end{lemma}
{\it Proof.} The group $G$ acts on a tree $T=(V,E)$ with chain complex
$$0\to\bigoplus_{\beta\in E/G}\C[G/G_\beta]\to\bigoplus_{\alpha\in V/G} \C [G/G_\alpha]\to \C.$$
Take a projective resolution of this complex in the category of 
chain complexes over $\cg$, say $P_*\to Q_* \to R_*$,
with $P_*$ a projective resolution for $\bigoplus \C[G/G_\beta]$, $Q_*$ one
for $\bigoplus \C[G/G_\alpha]$ and $R_*$ for $\C$. Upon
tensoring with $\vNG\otimes_{\cg}-$ we obtain
a short exact sequence of chain
complexes
$$\vNG\otimes_{\cg} P_*\to\vNG\otimes_{\cg} Q_*\to\vNG\otimes_{\cg} R_*;$$
the exactness results from the fact that the sequences
$P_i\to Q_i\to R_i$ are split exact for all $i$, because $R_i$ is projective.
Taking homology yields a long exact sequence of $L^{2}$-homology groups
$$\cdots\to H_{i+1}(G;\vNG)\to 
\bigoplus_{\beta\in E/G} H_{i}({\operatorname{Ind}}_{G_\beta}^G EG_\beta; {\vNG})\to$$
$$\bigoplus_{\alpha\in V/G} H_{i}({\operatorname{Ind}}_{G_\alpha}^G 
EG_\alpha;\vNG)\to H_{i}(G;\vNG)\to\cdots.$$
We used here that a $\cg$-projective resolution of $\C[G/G_\gamma]$ is 
chain homotopy equivalent to the cellular $\C G$-chain complex of the 
induced $G$-$CW$-complex 
$\operatorname{Ind}_{G_\gamma}^G EG_\gamma=G\times_{G_\gamma} EG_\gamma$. 
Therefore
$$H_*(\vNG\otimes_{\cg}P_*)= \bigoplus_{\beta\in E/G}
 H_*({\operatorname{Ind}}_{G_\beta}^G EG_\beta;\vNG)$$
and similarly for $H_*(\vNG\otimes_{\cg}Q_*)$.
According to \cite{L} (Theorem 6.54, (7)), 
for any induced $G$-$CW$-complex $\operatorname{Ind}_{G_\gamma}^G X$ 
one has
$$ \dvNG H_i(\operatorname{Ind}_{G_\gamma}^G X;\vNG)=
\operatorname{dim}_{G_\gamma} H_i(X;\mathcal{N}(G_\gamma))$$
and it follows that
$\dvNG H_i(\operatorname{Ind}_{G_\gamma}^G EG_\gamma;\vNG)=
\beta_i(G_\gamma).$
Thus, by taking the alternating sum of $L^{2}$-Betti numbers 
in the long exact homology sequence above,
the desired
formula follows.\hfill QED

\medskip
There are groups $G$ of type FP over $\C$
containing centralizers $C_G(s)$ which are not of type
FP over $\C$. Such examples have first been constructed 
by Leary and Nucinkis in \cite{LN},
and those cannot satisfy Brown's formula, because then $e(C_G(s))$
is not defined. 
The following group $G$ is a simple example for which Formula (\ref{main}) 
holds whereas (\ref{brown}) doesn't apply. 
Take first a group $\G$ as described by Leary-Nucinkis in \cite{LN} 
with the following property:

\medskip
\emph{$\G$ is of type FP over $\C$ and contains an element $t\in\G$ 
of finite order such that $C_{\G}(t)$ is not 
of type FP over $\C$.}

\medskip
\noindent
Then the right hand side of Brown's formula (\ref{brown})
doesn't make sense for the  
group 
$G= \G\times\Z$, which is of type FP over $\C$
but none of the centralizers $C_G((t,n))=C_\G(t)\times\Z$ are;
note that $(t,n)\in G$ is of finite order if and only if
$n= 0$.
But nevertheless, the group $G$ satisfies Conjecture \ref{main} because 
of the following.
\begin{lemma}\label{HtimesZ} Let $H$ be a group of type FP over $\C$ and
$G:=H\times\Z$. Then $G$ is of type FF over $\C$,
satisfies Conjecture \ref{conjecture} and
$W(G)=0\in K_0(\cg)$.\end{lemma}
{\it Proof.} Let $P_*:0\to P_n\to\dots\to P_1\to P_0\to\C\to 0$
be a resolutions of type FP over $\C$ for $H$ and $D_*:0\to\C\left<z\right>\to\C\left<z\right>\to \C\to 0$
be the projective resolution for $\Z=\left<z\right>$ with the map 
$\C\left<z\right>\to\C\left<z\right>$ given by 
multiplication with $1-z$. Then
$E_*=P_*\otimes D_*\to\C\to 0$ 
is a resolution of type FP over $\C$ for $G=H\times\Z$, and since
$E_i=(P_i\otimes\C\left<z\right>)\oplus (P_{i-1}\otimes\C\left<z\right>)$, 
we see that
$W(G)=\sum_{i=0}^{n+1}(-1)^i[E_i]=0$
(terms cancel pairwise); hence $G$ is of type FF
over $\C$ and $E(G)=0$ so that
$E(G)(s)=0$ for every $s\in G$. On the other hand, the
centralizer of $s=(u,v)\in H\times\Z$ contains the normal
subgroup $\{1_H\}\times\Z$ so that $\lte(C_G(s))=0$ as well.
\hfill QED.

\medskip
It follows that Conjecture \ref{conjecture} holds for any group of type FP over
$\C$ of the form $H\times\Z$, because both sides are zero; we shall 
construct non-zero examples later (recall that if $H\times \Z$
is of type FP over $\C$ then so is $H$ by Proposition 2.7
of \cite{Bieri}). We we will show in 
Section~\ref{Binfini} (Theorem \ref{rational}) that for each 
$\rho\in \Q$ 
there exists a group $G(\rho)$ of type 
FP over $\C$ containing an element $s$ of finite order such that 
$C_{G(\rho)}(s)$ 
is not of type FP over $\C$ but such that $G(\rho)$ nevertheless satisfies 
Conjecture \ref{conjecture}, with $E(G(\rho))(s)=\rho$.
\section{Stability properties of Formula (\ref{main})}\label{stability}
In this section we shall study some stability properties of Formula~(\ref{main}), starting with the following.
\begin{lemma}\label{products}
Let $A$ and $B$ be two groups of type FP over $\C$ such that $A$ satisfies
Formula~(\ref{main}) for some $a\in A$, and $B$ satisfies it for some $b\in B$.
Then $G=A\times B$ satisfies Formula~(\ref{main}) for the element $(a,b)$.\end{lemma}
{\it Proof.} Let $P_*:0\to P_n\to\dots\to P_1\to P_0\to\C\to 0$
be a resolution of type FP over $\C$ for $A$ and
$Q_*:0\to Q_n\to\dots\to Q_1\to Q_0\to\C\to 0$
one for $B$ (by adding trivial modules we can assume that both resolutions have the same length). A projective resolution of type FP over $\C$ for $G=A\times B$ is given by $E_*=P_*\otimes Q_*\to\C\to 0$. For an element $s=(a,b)\in G$ we compute
\begin{eqnarray*}HS(W(G))(s)&=&\sum_{i=0}^{2n}(-1)^iHS(E_i)(s)=\sum_{i=0}^{2n}(-1)^i\sum_{k+\ell=i}HS(P_k\otimes Q_{\ell})(a,b)\\
&=&\sum_{i=0}^{2n}\sum_{k+\ell=i}(-1)^{k+\ell}HS(P_k)(a)HS(Q_{\ell})(b)\\
&=&HS(W(A))(a)HS(W(B))(b)=\lte(C_A(a))\lte(C_B(b))\\
&=&\lte(C_A(a)\times C_B(b)).\end{eqnarray*}
Here we used in the second line the fact that for $P$ and $Q$ finitely generated projective modules over $\C A$ and $\C B$ respectively, 
$HS(P\otimes Q)(a,b)=HS(P)(a)HS(Q)(b).$
We conclude using Lemma \ref{prodchi} and
the fact that $C_G(a,b)$ = $C_A(a)\times C_B(b)$.
Note that the $L^2$-Betti numbers of $C_A(a)$ 
are finite, and trivial for large degrees, because $C_A(a)$ is 
assumed to have a well-defined $L^2$-Euler
characteristic and $\operatorname{cd}_\C C_A(a)$ is finite; a similar remark
applies to $C_B(b)$.\hfill QED
\begin{defi}[Condition (F)] The fundamental group $G$ of a 
finite graph of groups $\{G_\gamma\}$ satisfies \emph{Condition (F)}, 
if the $G$-action on the associated standard tree $T$ is such that
for every element of finite order $s\in G$, the action of $C_G(s)$ 
on the fixed tree $T^s$ satisfies the hypothesis
of Lemma \ref{graph}.
\end{defi}
\begin{rem}Condition (F) amounts to say that for any element 
of finite order $s\in G$ and for each of the stabilizers 
$H<C_G(s)$ appearing on the fixed tree $T^s$, 
the series $\sum_i\beta_i(H)$ is convergent and equals $0$ 
for all but finitely many conjugacy classes $(H)$.\end{rem}
\begin{lemma}\label{graphgroups} Let $G$ be the fundamental 
group of a finite graph of groups.
\begin{itemize}
\item[(i)]If all edge and vertex groups 
satisfy Formula (\ref{main}) at all elements of infinite order,
then so does $G$. 
\item[(ii)]If $G$ satisfies Condition (F) and
all edge and vertex groups satisfy Formula (\ref{main})
at all elements of finite order, 
then $G$ satisfies Formula (\ref{main})
at all elements of finite order.\end{itemize}\end{lemma}
{\it Proof.} The group $G$ acts cocompactly on a tree $T=(V,E)$, yielding a resolution 
$0\to\bigoplus_{\beta\in E/G}\C[G/G_{\beta}]\to
\bigoplus_{\alpha\in V/G}\C[G/G_{\alpha}]\to\C\to 0.$
Each of the groups $G_\gamma$ (for $\gamma\in V/G\sqcup E/G$) is of type FP 
over $\C$ (by assumption), so let us denote by 
$P^\gamma_*:0\to P^\gamma_n\to\dots\to P^\gamma_0\to\C\to 0$
a corresponding resolution of type FP. Tensoring by $\C[G/G_\gamma]$ 
yields the following resolutions of type FP over $\C$ of induced modules:
$\tilde{P}^\gamma_*: 0\to\tilde{P}^\gamma_n\to\dots\to
\tilde{P}^\gamma_1\to\tilde{P}^\gamma_0\to\C [G/G_{\gamma}]\to 0,\ 
{\text{for}}\,\, \gamma \in V/G\sqcup E/G,$
so that the Wall element for $G$ is given by
$$W(G)=\sum_{\alpha\in V/G}[\C [G/G_{\alpha}]]-\sum_{\beta\in E/G}[\C [G/G_{\beta}]],$$
where 
$[\C [G/G_{\gamma}]]=\sum_{j=0}^n(-1)^j[\tilde{P}^\gamma_j]
=i_*^\gamma W(G_\gamma)\in K_0(\cg).$
The complete Euler characteristic of $G$ is then given by
$$E(G)=\sum_{\alpha\in V/G}i_*^{\alpha}E(G_{\alpha})-
\sum_{\beta\in E/G}i_*^{\beta}E(G_{\beta}).$$
(i) Now let us take $s\in G$ of infinite order. By Cheeger-Gromov's Theorem~\ref{CG} of this note $\lte(C_G(s))=0$; on the other hand, $E(G)(s)=0$ 
because $E(G_\gamma)(t)=0$ for any $\gamma\in V/G\sqcup E/G$ and any $t$ of infinite order, by assumption on the $G_\gamma$'s.

\bigskip

\noindent
(ii) If $s\in G$ has finite order, then
\begin{eqnarray*}E(G)(s)&=&\sum_{[x]\in [s, G_{\alpha}]}E(G_{\alpha})(x)-
\sum_{[y]\in [s, G_{\beta}]}E(G_{\beta})(y)\\
&=&\sum_{[x]\in [s, G_\alpha]}\lte(C_{G_{\alpha}}(x))-
\sum_{[y]\in [s, G_{\beta}]}\lte(C_{G_{\beta}}(y))\end{eqnarray*}
because by assumption the $G_\gamma$'s satisfy Formula~(\ref{main}) at elements of finite order
(we used here the notation
$[s,G_\gamma]$ for the conjugacy classes of elements in $G_\gamma$
which are $G$-conjugate to $s$). 
So to conclude we need to show that the last line of the above equation is equal to $\lte(C_G(s))$, which we will do now.
We think of the $G_\gamma$'s as representatives for the 
stabilizers of the $G$-action on the standard tree $T$ of the given graph of groups
so that a general stabilizer will have the form
$t G_\gamma t^{-1}$.
Since $s$ has finite order, $T^s=(V^s,E^s)$ is a non-empty
tree, upon which $C_G(s)$ acts via the restriction of the $G$-action on $T$.
The stabilizer of a vertex or an edge $\in T^s$ has the form 
$C_G(s)\cap t G_\gamma t^{-1}$, where $s\in t G_\gamma t^{-1}$,
so that
$C_G(s)\cap t G_\gamma t^{-1} \cong C_{G_\gamma}(t^{-1}st).$
Moreover, by assumption $G$ satisfies Condition (F), 
and hence $\lte(C_{G_\gamma}(tst^{-1}))$ is well-defined so that $\lte(C_G(s))$ is
well defined too and, by Lemma \ref{graph} satisfies
$$\lte(C_G(s))=\sum_{x\in I}\lte(C_{G_\alpha}(x))-\sum_{y\in J}\lte(C_{G_\beta}(y))$$
with index set $I$ corresponding to $V^s/C_G(s)$. But this set corresponds 
bijectively to 
conjugacy classes of elements $x$ in the $[G_\alpha]$'s, 
which are $G$-conjugate to $s$; similarly for $J$.
\hfill QED
\section{Conjecture 1 and two classes of groups}\label{Binfini}
\smallskip
To begin with, we consider the following 
class $\mathcal{B}_\infty$
of groups.
\begin{defi} Let $\mathcal{B}_\infty$ denote the smallest class of groups which
contains all groups of type FF over $\C$, all groups of type FP over $\C$ 
which satisfy the Bass Conjecture over $\C$, 
all groups $G$ with cocompact $\underbar{E}G$,
all groups $G=H\times \Z$ with $H$ of type FP over $\C$
and which is closed
under finite products of groups and under passing to the fundamental
group of a finite graph of groups.
\end{defi}
Clearly all groups in $\mathcal{B}_\infty$ are of type FP over $\C$.
In particular, the Wall element $W(G)\in K_0(\C G)$ is defined for 
all groups in $\mathcal{B}_\infty$. Examples of groups in 
$\mathcal{B}_\infty$ include word hyperbolic groups, 
braid groups, cocompact CAT(0)-groups, Coxeter groups, 
mapping class groups of surfaces, knot groups, finitely 
generated one-relator groups, $S$-arithmetic groups, 
Artin groups, amenable groups of type FP over $\C$.
Many more groups can be obtained using the closure properties 
mentioned before; the groups thus obtained are in general 
not known to satisfy the Bass conjecture over $\C$. We do not 
know of any group of type FP over $\C$ not belonging 
to $\mathcal{B}_\infty$. As we have seen, 
there are groups $G$ in $\mathcal{B}_\infty$ containing $x$ of finite 
(resp. infinite) order, whose centralizer $C_G(x)$ is not of type FP 
over $\C$ and, therefore, $E(C_G(x))$ is not defined and 
$C_G(x)\not\in\mathcal{B}_\infty$. But nevertheless, the following holds.
\begin{thm}
Let $G$ be a group in $\mathcal{B}_\infty$ and $s\in G$ an element
of infinite order. Then Formula (\ref{main}) holds at $s$:
$$E(G)(s)=0=\chi^{(2)}(C_G(s)).$$
\end{thm}
{\it Proof.} We have already seen that the right 
hand side is $0$ (cf.~Cheeger-Gromov's 
Theorem~\ref{CG} in this note). The left hand side is certainly $0$ in 
case $G$ is of type FF over $\C$ or if $G$
satisfies the Bass Conjecture over $\C$ or if $\underbar{E}G$ is cocompact.
 Moreover, by Lemmas \ref{products} and \ref{graphgroups} (i), if $G=H\times K$ or $G$ is the fundamental group of a finite graph of groups $G_\alpha$ and if $E(L)(t)=0$ for all $t$ of infinite order in $L$, where $L$ is one of the groups $H, K, G_\alpha$, then $E(G)(s)=0$ for all elements of infinite order $s\in G$. Finally, $G=H\times\Z$ certainly satisfies $E(G)(s)=0$ for all $s$ (see Lemma \ref{HtimesZ}).\hfill QED
\smallskip

We now describe a class of groups $\mathcal{B}\subset\mathcal{B}_\infty$
containing many examples
of groups $G$ with $E(G)(s)\ne 0$ for some $s\ne e$ in $G$
satisfying Conjecture 
\ref{conjecture}, but such that the corresponding centralizer $C_G(s)$
is not of type FP over $\C$. 

\begin{defi} Let $\mathcal{B}$ denote the smallest class of groups 
which
contains all groups $G$ with cocompact $\underline{E}G$, 
all groups $G=H\times \Z$ with $H$ of type FP over $\C$ and which is closed
under finite products of groups 
and under passing to the fundamental
group of a finite graph of groups which satisfy Condition~(F).\end{defi}
\begin{thm}
The groups of the class $\mathcal{B}$ satisfy Conjecture \ref{conjecture}.
\end{thm}
{\it Proof.} This follows by applying Lemmas \ref{HtimesZ}, \ref{products}
and \ref{graphgroups}.\hfill QED
\begin{cor}\label{chis}
For $G$ satisfying Conjecture \ref{conjecture}, 
$\chi(G)=\sum_{[s]\in[G]}\lte(C_G(s)).$\end{cor}
{\it Proof.} By definition we have that 
$$\chi(G)=\sum_i(-1)^i\operatorname{dim}_\C H_i(G;\C)=
\sum_i(-1)^i\operatorname{dim}_\C \C\otimes_{\C G}P_i,$$
where
$P_*\to\C$ is a resolution of $G$ of type FP over
$\C$. It implies that $\sum_{[s]\in[G]} E(G)(s)=\chi(G)$, because
for $P$ a finitely generated projective
$\C G$-module, $\sum_{[s]\in[G]}\HS(P)(s)=\operatorname{dim}_\C
\C\otimes_{\cg}P$. The desired result now follows
from Formula (\ref{main}). 
\hfill QED

\medskip

\noindent
We shall now construct explicit non-trivial examples in 
the class ${\mathcal B}$. More precisely we prove the following.
\begin{thm}\label{rational} Given $\rho\in\Q$ there exists a group $G=G(\rho)$
of type FP over $\C$ with $s\in G$ of order $2$ such that $G$ satisfies Conjecture \ref{conjecture}, with
$$E(G)(s)=\lte(C_{G}(s))=\rho$$
but with the centralizer $C_{G}(s)$ not of type FP over $\C$.
\end{thm}
Before proceeding with the proof we need the following.
\begin{lemma}\label{rho} For $\rho\in\Q$ there exist a group $G_{\rho}\in\mathcal{B}$ with $\lte(G_{\rho})=\rho$.\end{lemma}
{\it Proof.} Since a free group $F_n$ of rank $n$ satisfies $\lte(F_n)=1-n$, one has for $n,k\ge 0$ that
$\lte((F_2\times F_{n+1})\ast F_k)=n-k$,
so that for $\ell>0$
$$\lte(((F_2\times F_{n+1})\ast F_k)\times \Z/\ell\Z)=\frac{n-k}{\ell}.$$
The group $G=((F_2\times F_{n+1})\ast F_k)\times\Z/\ell\Z$ admits a cocompact $\underline{E}G$ via its obvious quotient action on $E(G/(\Z/\ell\Z))$, with orbit space the finite complex $((\vee^2 S^1)\times(\vee^{n+1}S^1))\vee (\vee^kS^1),$ thus $G\in\mathcal{B}$.\hfill QED

\medskip
{\it Proof of Theorem \ref{rational}.} Let $\G$ be one of the groups described in \cite{LN}, Example 9, such
that $\G$ is of type FP over $\C$, $s\in\G$ is an element of order $2$ and $C_{\G}(s)$ is not finitely generated. By definition of $\mathcal{B}$, the group $H:=\G\times\Z$ belongs to $\mathcal{B}$,
and $C_H((s,0))$ is not finitely generated, because it maps onto $C_{\G}(s)$. Writing $t$ for $(s,0)$, we form
$K:=H\ast_{\left<t\right>}H\in\mathcal{B}.$
Thus, $K$ is the fundamental group of a finite graph of groups $\{H,\left<t\right>\}$, with associated tree $T$. If $w\in K$ has finite order with $w$ not conjugate to $t$, the edge stabilizers of the $C_K(w)$ action on $T^w$ are all trivial, and the vertex stabilizers are isomorphic to $C_H(z)$ for some element $z$ of order 2 in $H$, thus $\beta_i(C_H(z))=0$ for all $i$,
because such a centralizer contains a normal infinite cyclic subgroup. The centralizer of $\left<t\right>$ in $K$ decomposes as a fundamental group of a graph of groups of the form $\{H_\delta,\left<t\right>\}$ with the $H_\delta$'s again isomorphic to groups $C_H(w)$, $w\in H$
of order 2, so that $\beta_i(H_\delta)=0$ for all $i$ and all $\delta$. It follows that $K$ satisfies Condition (F) and
$\lte(C_K(t))=-\lte(\left<t\right>)= -{1}/{2}.$
Note that $C_K(t)/{\left<t\right>}$ maps onto $C_H(t)/{\left<t\right>}$, which shows
that $C_K(t)$ is not finitely generated.
Forming
$G:= K\times G_{-2\rho}\in\mathcal{B}$
where $G_{-2\rho}$ is obtained following Lemma \ref{rho} above,
gives a group with $C_{G}(t)=C_K(t)\times G_{-2\rho}$ not of type
FP over $\C$ (because it is not finitely generated), but
$$\lte(C_{G}(t))=\lte(C_K(t))\cdot\lte(G_{-2\rho})=-\frac{1}{2}\cdot (-2\rho)=\rho.$$
\hfill QED

\bigskip

{\it Acknowledgements.} The first named author thanks Ken Brown for friendly discussions on \cite{Brown}.

\end{document}